\documentclass[12pt]{article}
\usepackage{amssymb, amsmath, amsthm, url, graphicx}

\def\3{\subset }
\def\4{\subseteq }
\def\<{\left<}
\def\>{\right>}

\def\bit{\begin{itemize}}
\def\eit{\end{itemize}}
\def\3{\subset }
\def\4{\subseteq }
\def\ov{\overline}

\def\0{\leqno}

\def\barr{\begin{array}}
\def\earr{\end{array}}

\def\Z{{\rlap{$\kern2pt{\rm Z}$}{\rm Z}\,}}

\title{\bf Finite groups determined\\ by an inequality of the orders\\ of their elements}
\author{Marius T\u arn\u auceanu}
\date{October 1, 2012}

\begin{document}

\maketitle

\begin{abstract}
In this note we introduce and characterize a class of finite
groups for which the element orders satisfy a certain inequality.
This is contained in some well-known classes of finite groups.
\end{abstract}



\section{Introduction}

Let ${\rm CP}_1$, CP and CN be the classes of finite groups in
which the centralizers of all nontrivial elements contain only
elements of prime order, of prime power order and are nilpotent,
respectively. Clearly, we have ${\rm CP}_1 \subset {\rm CP}
\subset {\rm CN}$. Moreover, the classes ${\rm CP}_1$ and CP
consist of exactly those finite groups all of whose elements have
prime order and prime power order, respectively. They have been
studied in many papers, as \cite{1}-\cite{3}, \cite{5}-\cite{7}
and \cite{13}.

In the following we consider the finite groups $G$ such that
$$o(xy) \leq {\rm max}\{o(x), o(y)\}, \mbox{ for all } x, y \in
G. \0(*)$$These form another interesting subclass of CP, that will
be denoted by ${\rm CP}_2$. Its exhaustive description is the main
goal of this note.

Most of our notation is standard and will not be repeated here.
Basic notions and results on group theory can be found in
\cite{4,8,9,11}.

First of all, we observe that if a finite group $G$ belongs to
${\rm CP}_2$, then for every $x, y \in G$ satisfying $o(x) \neq
o(y)$ we have $$o(xy)={\rm max}\{o(x), o(y)\},$$that is the order
map is very close to a monoid homomorphism from $(G, \cdot)$ to
$(\mathbb{N}^*, {\rm max})$.
\bigskip

An immediate characterization of finite groups contained in ${\rm
CP}_2$ is indicated in the following theorem.

\bigskip\noindent{\bf Theorem A.} {\it Let $G$ be a finite group and set $\pi_e(G)=\{o(x) \mid x \in G \}$.
Then the following conditions are equivalent:
\begin{itemize}
\item[\rm a)]
        $G$ belongs to ${\rm CP}_2$.
\item[\rm b)]
        For every $\alpha \in \pi_e(G)$, the set $G_{\alpha}=\{x \in G \mid o(x) \leq \alpha\}$ is a normal subgroup of $G$.
\end{itemize}}

Next, we will focus on establishing some connections between ${\rm
CP}_2$ and the previous classes CP and ${\rm CP}_1$.

\bigskip\noindent{\bf Proposition B.} {\it The class ${\rm CP}_2$ is properly contained in the class ${\rm CP}$.}
\bigskip

On the other hand, by taking $\sigma = (12)(34), \tau = (235) \in
A_5$, one obtains $$5=o(\sigma\tau) > 3 = {\rm max}\{o(\sigma),
o(\tau)\},$$and therefore ${\rm CP}_2$ does not contain the
alternating group $A_5$. Since $A_5$ belongs to ${\rm CP}_1$, we
conclude that ${\rm CP}_1$ is not contained in ${\rm CP}_2$. It is
obvious that the converse inclusion also fails (for example, any
abelian $p$-group belongs to ${\rm CP}_2$, but not to ${\rm
CP}_1$).

\bigskip\noindent{\bf Remarks.}
\begin{itemize}
\item[\rm 1.]
        Other two remarkable classes of finite $p$-groups, more large as
the class of abelian $p$-groups, are contained in ${\rm CP}_2$:
regular $p$-groups (see Theorem 3.14 of \cite{11}, II, page 47)
and $p$-groups whose subgroup lattices are modular (see Lemma
2.3.5 of \cite{10}). Moreover, by the main theorem of \cite{12},
we infer that the powerful $p$-groups for $p$ odd also belong to
${\rm CP}_2$.
\item[\rm 2.]
        The smallest nonabelian $p$-group contained in ${\rm
CP}_2$ is the quaternion group $Q_8$, while the smallest $p$-group
not contained in ${\rm CP}_2$ is the dihedral group $D_8$. Notice
that all quaternion groups $Q_{2^n}$, for $n \geq 4$, as well as
all dihedral groups $D_n$, for $n \neq 1, 2, 4$, are not contained
in ${\rm CP}_2$.
\item[\rm 3.]
        The class ${\rm CP}_2$ contains finite groups which are
not $p$-groups, too. The smallest example of such a group is
$A_4$. Remark that the groups $A_n$, $n \geq 5$, does not belong
to ${\rm CP}_2$, and this is also valid for the symmetric groups
$S_n$, $n \geq 3$.
\end{itemize}

Clearly, ${\rm CP}_2$ is closed under subgroups. On the other
hand, the above results imply that ${\rm CP}_2$ is not closed
under direct products or extensions. The same thing can be said
with respect to homomorphic images, as shows the following
example.

\bigskip\noindent{\bf Example.} Let $p$ be a prime and $G$ be the
semidirect product of an elementary abelian $p$-group $A$ of order
$p^p$ by a cyclic group of order $p^2$, generated by an element
$x$ which permutes the elements of a basis of $A$ cyclically. Then
it is easy to see that $G$ belongs to ${\rm CP}_2$, $x^p\in Z(G)$
and the quotient $Q=\frac{G}{\langle x^p \rangle}$ is isomorphic
to a Sylow $p$-subgroup of $S_{p^2}$. Obviously, in $Q$ a product
of two elements of order $p$ can have order $p^2$, and hence it
does not belong to ${\rm CP}_2$.
\bigskip

The next result collects other basic properties of the finite
groups contained in ${\rm CP}_2$.

\bigskip\noindent{\bf Proposition C.} {\it Let $G$ be a finite group contained in ${\rm CP}_2$. Then:
\begin{itemize}
\item[\rm a)]
        There is a prime $p$ dividing the order of $G$ such that $F(G)=O_p(G)$.
\item[\rm b)]
        Both $Z(G)$ and $\Phi(G)$ are $p$-groups.
\item[\rm c)]
        $Z(G)$ is trivial if $G$ is not a $p$-group.
\end{itemize}}

We are now able to present our main result, that gives a complete
des\-crip\-tion of the class ${\rm CP}_2$.

\bigskip\noindent{\bf Theorem D.} {\it A finite group $G$ is contained in ${\rm CP}_2$
if and only if one of the following statements holds:
\begin{itemize}
\item[\rm a)]
        $G$ is a $p$-group and $\Omega_n(G) = \{x \in G \mid x^{p^{n}}=1\}$, for all $n \in \mathbb{N}$.
\item[\rm b)]
        $G$ is a Frobenius group of order $p^{\alpha}q^{\beta}$,
        $p < q$, with kernel $F(G)$ of order $p^{\alpha}$ and
        cyclic complement.
\end{itemize}}

Since all $p$-group and all groups of order $p^{\alpha}q^{\beta}$
are solvable, Theorem D leads to the following corollary.

\bigskip\noindent{\bf Corollary E.} {\it The class ${\rm CP}_2$ is properly contained in the class of finite sol\-va\-ble groups.}

\bigskip\noindent{\bf Remark.} The finite supersolvable groups and
the CLT-groups constitute two important subclasses of the finite
solvable groups. Since $A_4$ belongs to ${\rm CP}_2$, we infer
that ${\rm CP}_2$ is not included in these classes. Conversely, a
finite supersolvable group or a CLT-group does not necessarily
possess the structure described above, and thus they are not
necessarily contained in the class ${\rm CP}_2$.
\bigskip

As we already have seen, both ${\rm CP}_1$ and ${\rm CP}_2$ are
subclasses of CP, and each of them is not contained in the other.
Consequently, an interesting problem is to find the intersection
of these subclasses. This can be made by using again Theorem D.

\bigskip\noindent{\bf Corollary F.} {\it A finite group $G$ is contained
in the intersection of \,${\rm CP}_1$ and ${\rm CP}_2$ if and only
if one of the following statements holds:
\begin{itemize}
\item[\rm a)]
        $G$ is a $p$-group of exponent $p$.
\item[\rm b)]
        $G$ is a Frobenius group of order $p^{\alpha}q$,
        $p < q$, with kernel $F(G)$ of order $p^{\alpha}$ and exponent $p$,
        and cyclic complement.
        Moreover, in this case we have $G'=F(G)$.
\end{itemize}}

\medskip\noindent{\bf Remark.} $A_4$ is an example of a group
of type b) in the above corollary. Mention that for such a group
$G$ the number of Sylow $q$-subgroups is $p^{\alpha}$. It is also
clear that $G$ possesses a nontrivial partition consisting of
Sylow subgroups: $F(G)$ and all conjugates of a Frobenius
complement.
\bigskip

Finally, we indicate a natural problem concerning the class of
finite groups introduced in our paper.

\bigskip\noindent{\bf Open problem.} Give a precise description of
the structure of finite $p$-groups contained in ${\rm CP}_2$.

\section{Proofs of the main results}

\bigskip\noindent{\bf Proof of Theorem A.} Assume first that $G$
belongs to ${\rm CP}_2$. Let $\alpha \in \pi_e(G)$ and $x, y \in
G_{\alpha}$. Then, by $(*)$, we have $$o(xy) \leq {\rm max}\{o(x),
o(y)\} \leq \alpha,$$which shows that $xy \in G_{\alpha}$. This
proves that $G_{\alpha}$ is a subgroup of $G$. Moreover,
$G_{\alpha}$ is normal in $G$ because the order map is constant on
each conjugacy class.

Conversely, let $x, y \in G$ and put $\alpha = o(x), \beta =
o(y)$. By supposing that $\alpha \leq \beta$, one obtains $x, y
\in G_{\beta}$. Since $G_{\beta}$ is a subgroup of $G$, it follows
that $xy \in G_{\beta}$. Therefore $$o(xy) \leq \beta = {\rm
max}\{o(x), o(y)\},$$completing the proof.
\hfill\rule{1,5mm}{1,5mm}

\bigskip\noindent{\bf Proof of Proposition B.} Let $G$ be finite
group in ${\rm CP}_2$ and take $x\in G$. It is well-known that $x$
can be written as a product of (commuting) elements of prime power
orders, say $x=x_1x_2\cdots x_k$. Then the condition $(*)$ implies
that $$\prod_{i=1}^k o(x_i)=o(x)\leq {\rm max}\{o(x_i)\mid
\,i=\ov{1,k}\,\},$$and so $k=1$. Hence $x$ is of prime power
order, i.e. $G$ is contained in CP.

Obviously, the inclusion of ${\rm CP}_2$ in CP is strict (we
already have seen that $A_5$ belongs to CP, but not to ${\rm
CP}_2$). \hfill\rule{1,5mm}{1,5mm}

\bigskip\noindent{\bf Proof of Proposition C.}
\begin{itemize}
\item[\rm a)] We know that $F(G)$ is the product of the
subgroups $O_p(G)$, where $p$ runs over the prime divisors of
$\mid G \mid$. Suppose that there are two distinct primes $p$ and
$q$ dividing the order of $F(G)$. This leads to the existence of
two elements $x$ and $y$ of $F(G)$ such that $o(x) = p$ and $o(y)
= q$. Since $F(G)$ is nilpotent, we obtain $xy = yx$ and so $o(xy)
= pq$, a contradiction. Thus $F(G) = O_p(G)$, for a prime divisor
$p$ of $\mid G \mid$.
\item[\rm b)] It is well-known that both $Z(G)$ and $\Phi(G)$ are
normal nilpotent subgroup of $G$. By the maximality of $F(G)$, it
follows that $Z(G)$ and $\Phi(G)$ are contained in $F(G)$, and
therefore they are also $p$-groups.
\item[\rm c)] Assume that $Z(G)$ is not trivial and take $x\in Z(G)$
with $o(x)=p$. If $G$ is not a $p$-group, it contains an element
$y$ of prime order $q\neq p$. Then $o(xy)=pq$, contradicting
Proposition B.
\hfill\rule{1,5mm}{1,5mm}
\end{itemize}

\bigskip\noindent{\bf Proof of Theorem D.} If $G$ is a $p$-group,
then the conclusion is obvious.\

Assume now that $G$ is not a $p$-group. We will proceed by
induction on $|G|$. Since $G$ belongs to ${\rm CP}_2$, all the
numbers in $\pi_e(G)$ are prime powers. Let $q^n$ be the largest
number of $\pi_e(G)$, where $q$ is a prime, and let $N=\{g\in
G\mid o(g)<q^n\}$. Then $N\unlhd G$ and ${\rm exp}(G/N)=q$. Since
$|N|<|G|$, by the inductive hypothesis it follows that either $N$
is a $p$-group or $N$ is a Frobenius group with kernel $K$ of
order $p^{\alpha}$ and cyclic complement $H$ of order $r^{\beta}$,
where $p, r$ are distinct primes. We will prove that in both cases
$G$ is a Frobenius group whose kernel and complement are
$p$-groups.

\bigskip\hspace{10mm}{\bf Case 1.} $N$ is a $p$-group.

\noindent Since $G$ is not a $p$-group, we can take $Q\in
Syl_q(G)$, where $p\neq q$. So $G=N\rtimes Q$. Since every element
of $N$ is of prime power order, we have $C_N(h)=1$ for all $1\neq
h\in Q$. Thus, $G$ is a Frobenius group with kernel $N$ and
complement $Q$.

\bigskip\hspace{10mm}{\bf Case 2.} $N$ is a Frobenius group.

\bigskip\hspace{15mm}{\bf Subcase 2.1.} $q\neq p$ and $q\neq r$.

\noindent By a similar argument as that of Case 1, we know that
$G$ is a Frobenius group with kernel $N$. But $N$ is not
nilpotent, a contradiction.

\bigskip\hspace{15mm}{\bf Subcase 2.2.} $q=r$.

\noindent Let $Q\in Syl_q(G)$. Then $G=K\rtimes Q$. By a similar
argument as that of Case 1, we know that $G$ is a Frobenius group
with kernel $K$ and complement $Q$.

\bigskip\hspace{15mm}{\bf Subcase 2.3.} $q=p$.

\noindent We observe that all elements of $G\setminus N$ are of
order $q^n$, and $g^q\in K$, where $g\in G\setminus N$ and $K\in
Syl_q(N)$. So if $N_G(H)\cap(G\setminus N)\neq 1$, then
$N_G(H)\cap K\neq 1$. But $N$ is a Frobenius group and $N_N(H)=H$.
It follows that $N_G(H)=H$. Since $H$ is cyclic, $N_G(H)=C_G(H)$.
One obtains that $G$ is $r$-nilpotent and thus $G=P\rtimes H$,
where $P\in Syl_p(G)$. A similar argument as that of Case 1 shows
that $G$ is again a Frobenius group with kernel $P$ and complement
$H$.

Finally, we prove that $H$ is cyclic. By Burnside's Theorem we
only need to prove that $H$ is not a 2-group. If not, let
$L=\{g\in G\mid o(g)=2\}$. Then $L\unlhd G$. It follows that
$K\times L\leq G$, where $K$ is the Frobenius kernel. This
contradicts the fact that all elements of $G$ are of prime power
order. \hfill\rule{1,5mm}{1,5mm}

\bigskip\noindent{\bf Proof of Corollary F.} The equivalence follows
directly by Theorem D. In this way, we have to prove only that
$G'=F(G)$ in the case b).

Obviously, $G' \subseteq F(G)$. For the converse inclusion, let $x
\in F(G)$ be a nontrivial element. Then $o(x) = p$. If $y$ is an
arbitrary element of order $q$ in $G$, then we have $$o(xy) \leq
{\rm max}\{o(x), o(y)\} = q,$$and therefore $o(xy) \in \{p, q\}$.
If we assume that $o(xy) = p$, it results
$$q = o(y) = o(x^{-1}xy) \leq {\rm max} \{o(x^{-1}), o(xy)\} =
p,$$a contradiction. This shows that $o(xy) = q$. Then there is $z
\in G$ such that $xy \in \langle y\rangle^z$, say $xy =
z^{-1}y^kz$ with $k \in \mathbb{Z}$. Since the element
$$xy^{1-k}=z^{-1}y^kzy^{-k}=[z,y^k]$$has order $p$, we infer that
$k$ must be equal to 1. Hence $$x = [z,y] \in G',$$which completes
the proof.
\hfill\rule{1,5mm}{1,5mm}
\bigskip
\bigskip

{\bf Acknowledgements.} The author is grateful to the reviewer for
its remarks which improve the previous version of the paper.

\vspace*{4ex}\small

\hfill
\begin{minipage}[t]{5cm}
Marius T\u arn\u auceanu \\
Faculty of  Mathematics \\
``Al.I. Cuza'' University \\
Ia\c si, Romania \\
e-mail: {\tt tarnauc@uaic.ro}
\end{minipage}

\end{document}